\documentclass[12pt,a4paper]{article}
\usepackage{amssymb}
\usepackage{amsmath}
\usepackage{color}
\usepackage[colorlinks,linkcolor=blue,citecolor=red]{hyperref}

\newcommand{\Hess}{\operatorname{Hess}}
\newcommand{\bbbr}{\mathbb R}
\newcommand{\bbbc}{\mathbb C}

\newtheorem{theorem}{Theorem}

\begin{document}

\title{On the integrability of symplectic Monge-Amp\'ere equations}

\author{B. Doubrov and E.V. Ferapontov}
     \date{}
     \maketitle
     \vspace{-5mm}
\begin{center}
Department of Mathematical Physics\\
Faculty of Applied Mathematics\\
 Belarussian State University\\
  Nezavisimosti av. 4, 220030 Minsk, Belarus\\
  \ \\
  and\\
  \ \\
Department of Mathematical Sciences \\ Loughborough University \\
Loughborough, Leicestershire LE11 3TU \\ United Kingdom \\[2ex]
e-mails: \\[1ex] \texttt{doubrov@islc.org} \\
\texttt{E.V.Ferapontov@lboro.ac.uk}\\

\end{center}

\vspace{1cm}
\begin{abstract}
Let $u$ be a function of $n$ independent variables $x^1, ..., x^n$, and
$U=(u_{ij})$  the Hessian matrix  of $u$. The symplectic Monge-Amp\'ere
equation is defined as a linear relation among all possible minors of  $U$.
Particular examples include the equation
%${\rm det} \ U=\triangle u$ arising in the theory of special Lagrangian three-folds in $\bbbc^n$, the equation
$\det U=1$ governing improper affine spheres and the so-called
heavenly equation, $u_{13}u_{24}-u_{23}u_{14}=1$, describing
self-dual Ricci-flat $4$-manifolds. In this paper we classify
integrable symplectic Monge-Amp\'ere equations in four dimensions
(for  $n=3$ the integrability of such equations   is known to be
equivalent to their  linearisability). This problem can be
reformulated geometrically as the classification of  `maximally
singular' hyperplane sections of the Pl\"ucker embedding of the
Lagrangian Grassmannian. We formulate a conjecture that any
integrable equation  of the form $F(u_{ij})=0$ in more than three
dimensions is necessarily of the symplectic Monge-Amp\'ere type.

\bigskip

\end{abstract}

%\vspace{-7mm}
\newpage

\section{Introduction}

Let us consider a function $u(x^1, ... x^n)$ of $n$ independent variables and
introduce the $n\times n$ Hessian matrix   $U=(u_{ij})$ of its second order
partial derivatives. The symplectic Monge-Amp\'ere equation is a PDE of the
form
\begin{equation}
M_n+M_{n-1}+... +M_1+M_0=0 \label{M}
\end{equation}
where $M_l$ is a constant-coefficient linear combination of all
$l\times l$ minors of  the matrix $U$, $0\leq l\leq n$. Thus,
$M_n= \det U= \Hess u$, $M_0$ is a constant, etc. Equivalently,
these PDEs can be obtained by equating to zero a
constant-coefficient  $n$-form in the $2n$ variables $x^i, u_i$.
Equations of this type belong to the class of  completely
exceptional Monge-Amp\'ere equations  introduced in
\cite{Boillat}. Geometric and algebraic aspects of symplectic
Monge-Amp\'ere equations have been thoroughly investigated in
\cite{Lychagin, Banos}. We point out that the left hand side of
(\ref{M}), $M(U)$, constitutes the general form of null Lagrangian
densities, that is, functionals of the form $\int M(U) d{\bf x}$
which generate trivial Euler-Lagrange equations \cite{Olver}. The
class of equations (\ref{M}) is invariant under the natural
contact action of the symplectic group $Sp(2n)$, which is thus the
equivalence group of our problem. All subsequent classification
results will be formulated modulo this $Sp(2n)$-equivalence.

In the case $n=2$ one arrives at the standard Monge-Amp\'ere equations,
\begin{equation}
\epsilon(u_{11}u_{22}-u^2_{12})+\alpha u_{11}+\beta u_{12}+\gamma u_{22}+\delta
=0; \label{MAA}
\end{equation}
these are known to be the only equations of the form $F(u_{11},
u_{12}, u_{22})=0$ which are linearisable by a transformation from
the equivalence group $Sp(4)$.

The case $n=3$ is also understood completely: one can show that any
non-degenerate symplectic Monge-Amp\'ere equation is either linearizable, or
$Sp(6)$-equivalent to either of the canonical forms,
\begin{equation}
\Hess u=1,\qquad  \Hess u=u_{11}+u_{22}+u_{33}, \qquad  \Hess
u=u_{11}+u_{22}-u_{33},
 \label{MA}
\end{equation}
see  \cite{Lychagin, Banos}; we point out that all three canonical
forms are $Sp(6)$-equivalent over $\bbbc$. The first equation
arises in the theory of improper affine spheres, while the second
describes special Lagrangian $3$-folds in $\bbbc^3$ \cite{Calabi,
Joyce}. Here the non-degeneracy is understood as follows: let
$F(u_{ij})=0$ be a symplectic Monge-Amp\'ere equation. Consider
the linearized equation, $\partial F/\partial u_{ij} \ v_{ij}=0$,
obtained by setting $u\to u+\epsilon v$ and keeping  terms of the
order $\epsilon$. The non-degeneracy means that the corresponding
symbol, $\partial F/\partial u_{ij} \ \xi_i\xi_j$,  defines an
irreducible quadratic form.

The problem of  integrability of symplectic Monge-Amp\'ere equations was
addressed in \cite{Fer4} based on the method of hydrodynamic reductions
\cite{Fer1, Fer2,  Fer3}. Without going into technical details of this method,
let us formulate the main result needed for our purposes:
\begin{theorem} \label{thm1} \cite{Fer4}.
A non-degenerate three-dimensional symplectic Monge-Amp\'ere equation is
integrable by the method of hydrodynamic reductions if and only if it is
linearisable.
\end{theorem}
In particular, the  PDEs (\ref{MA}) are {\it not} integrable.
Although this result is essentially negative, it  will be crucial
for the classification of integrable equations in higher
dimensions (where the situation is far more interesting). Before
we proceed to the description of the main results, let us clarify
the geometry behind symplectic Monge-Amp\'ere equations (\ref{M})
and the linearisability/integrability conditions.  Let us consider
the Lagrangian Grassmannian $\Lambda$, which can be  (locally)
identified  with the space of $n\times n$ symmetric matrices $U$.
The minors of $U$ define the Pl\"ucker embedding of $\Lambda$ into
projective space $P^N$ (we identify $\Lambda$ with the image of
its projective embedding). Thus,  symplectic Monge-Amp\'ere
equations  correspond to  hyperplane sections of $\Lambda$. For
$n=3$ we have $\Lambda^6\subset P^{13}$, and linearisable
equations correspond to hyperplanes which are tangential to
$\Lambda^6$.
 Therefore, for $n=3$ the linearisability
condition  coincides with the equation of the dual variety of
$\Lambda^6$, which is known to be a hypersurface of degree four in
$(P^{13})^*$ (we refer to  \cite{Landsberg1, Landsberg2, Mukai}
for a general theory behind this example). In Sect. 2 we provide
the following characterisation of linearisable equations in any
dimension:

 \begin{theorem}\label{thm2} For a non-degenerate symplectic Monge-Amp\'ere equation (\ref{M})  the following conditions are equivalent:

 \noindent (1) The equation is linearisable by a transformation from $Sp(2n)$.

 \noindent (2) The equation is invariant under  an $n^2$-dimensional subalgebra of $Sp(2n)$.

 \noindent (3) The equation corresponds to a hyperplane which contains an osculating
 subspace  $O_{n-2}$  of the Lagrangian Grassmannian of the order $n-2$.
 \end{theorem}
%\begin{theorem}
%A symplectic Monge-Amp\'ere equation (\ref{M}) is linearisable by a transformation from
%the equivalence group $Sp(2n)$  if and only if the corresponding hyperplane contains
%an osculating subspace of the Lagrangian Grassmannian of the order $n-2$.
%\end{theorem}
For $n=3$ the third condition  reduces to the requirement of tangency.

In Sect. 3 we address the problem of  integrability of symplectic
Monge-Amp\'ere equations in four dimensions, $n=4$.
 Among the best known four-dimensional integrable examples one should
primarily mention the `heavenly' equation,
\begin{equation}
u_{13}u_{24}-u_{23}u_{14}=1, \label{heavenly}
\end{equation}
\cite{Plebanski} which is descriptive of Ricci-flat self-dual 4-manifolds
\cite{Atiyah}. It was demonstrated in \cite{Fer2} that this equation is
integrable by the method of hydrodynamic reductions. Although, in principle,
the method of hydrodynamic reductions can be applied in any dimension,  it
leads to a quite complicated analysis.  One way to bypass lengthy calculations
is based on the following simple idea: suppose our aim is the classification of
four-dimensional integrable equations of the form (\ref{M}) for a function
$u(x^1, x^2, x^3, x^4)$. Let us look for traveling wave solutions in the form
$$
u=u(x^1+\alpha x^4, \ x^2+\beta x^4, \ x^3+\gamma x^4),
$$
or, more generally,
$$
u=u(x^1+\alpha x^4, \ x^2+\beta x^4, \ x^3+\gamma x^4)+Q(x, x),
$$
where $Q$ is an arbitrary quadratic form in the variables $x^1, x^2, x^3, x^4$. The substitution of this ansatz into (\ref{M}) leads to a three-dimensional
symplectic Monge-Amp\'ere equation which must be integrable for {\it any}
values of constants $\alpha, \beta, \gamma$, and an arbitrary quadratic form $Q$. Since,  in  three dimensions, the
integrability conditions are explicitly known (and are equivalent to the
linearisability), this provides strong restrictions on the original
four-dimensional equation which are therefore {\it necessary} for the
integrability. In fact, in the present  context  they turn out to be
sufficient: if all three-dimensional equations obtained from a given
four-dimensional PDE by
traveling wave reductions are linearisable, then the PDE is integrable. The
philosophy of this approach is well familiar from the soliton theory: symmetry
reductions of integrable systems must be themselves integrable. Thus, for the
first heavenly equation, traveling wave solutions are governed by
$$
\alpha(u_{12}u_{13}-u_{11}u_{23})+\beta(u_{13}u_{22}-u_{12}u_{23})=1,
$$
which is a three-dimensional symplectic Monge-Amp\'ere equation.
One can show that it is indeed  linearisable for any values of
constants. This approach suggests a simple geometric
characterisation of integrable equations in four dimensions. Let
us first point out that for $n=4$ the Lagrangian Grassmannian
$\Lambda^{10}\subset P^{41}$ is foliated by a 7-parameter family
of $\Lambda^6$ where each $\Lambda^6$ corresponds to a collection
of Lagrangian planes in the symplectic space $V^8$ which pass
through a fixed vector.
\begin{theorem} \label{thm3}
A non-degenerate four-dimensional symplectic Monge-Amp\'ere equation is
integrable if and only if the corresponding hyperplane  is tangential to the
Lagrangian Grassmannian $\Lambda^{10}$ along a 4-dimensional subvariety which
meets all $\Lambda^6\subset \Lambda^{10}$.
\end{theorem}
Thus, integrable equations correspond to maximally singular hyperplane
sections. The classification of such hyperplanes leads to a complete list of integrable  equations in four
dimensions:
\begin{theorem}\label{thm4}
Over the field of complex numbers any integrable non-degenerate symplectic
Monge-Amp\'ere equation is $Sp(8)$-equivalent to one of the following  normal
forms:
\begin{enumerate}
\item $u_{11}-u_{22}-u_{33}-u_{44}=0$ (linear wave equation);
\item $u_{13}+u_{24}+u_{11}u_{22}-u_{12}^2=0$ (second heavenly equation);
\item $u_{13} = u_{12}u_{44}-u_{14}u_{24}$ (modified heavenly equation);
\item $u_{13}u_{24}-u_{14}u_{23}=1$ (first heavenly equation);
\item $u_{11}+u_{22}+u_{13}u_{24}-u_{14}u_{23}=0$ (Husain equation);
\item $\alpha u_{12}u_{34} + \beta u_{13}u_{24}+\gamma u_{14}u_{23} = 0$ (general
heavenly equation), $\alpha+\beta+\gamma=0$.
\end{enumerate}
\end{theorem}
Eqs. 2-6 are non-linearisable, and contact non-equivalent. They  have appeared  in
different contexts in  \cite{Plebanski, Husain, Lychagin}. To the best of our
knowledge, the integrability of Eqs. 3 and 6 was not recorded before. We
refer to \cite{Takasaki1, Takasaki2, Strachan1,  Grant, Dun1, Dun2, Nutku,
Sheftel,  Bogdanov1, Bogdanov2} for the Hamiltonian, twistorial  and symmetry aspects of
heavenly-type equations.

Symplectic Monge--Amp\'ere equations play a distiguished role in the classification
of all integrable PDEs of the `dispersionless Hirota'  type,
\begin{equation}
F(u_{ij})=0, \label{Hirota}
\end{equation}
which constitute a single relation among second order partial derivatives of a
function $u(x^1, ..., x^n)$. For $n=3$ integrable equations of this type were
studied  in \cite{Fer4, Odesskii2, Smith}, revealing a remarkable
correspondence with the theory of  hypersurfaces of the Lagrangian
Grassmannian,  generalised hypergeometric functions and $GL(2,
\bbbr)$-structures. Among them, only linearisable equations belong to the class
(\ref{M}). On the contrary, for $n\geq 4$ all known integrable examples belong
to the class (\ref{M}). This suggest the following conjecture:

\medskip

\noindent {\bf Conjecture.} {\it  For $n\geq 4$ any non-degenerate integrable
equation of the form (\ref{Hirota}) is necessarily of  the symplectic
Monge-Amp\'ere type.}

\medskip

\noindent Provided this conjecture is true, Theorem \ref{thm4} gives a complete
list of integrable equations of the form (\ref{Hirota}) in four dimensions.

\section{Characterisation of linearisable equations: proof of Theorem \ref{thm2}}

The results of this section are valid for  an arbitrary number of independent
variables   $n$. We begin with a few general remarks. The Lagrangian
Grassmannian $\Lambda$  is a variety of dimension $n(n+1)/2$ which can be
locally identified with $n\times n$ symmetric matrices $U=(u_{ij})$ so that
$u_{ij}$ can be viewed as local coordinates on $\Lambda$.  The action of
$Sp(2n)$  is defined by the infinitesimal operators
 $$
 X_{ij}=\frac{\partial}{\partial u_{ij}},
 $$
  $$
 L_{ij}=\sum_s u_{js}\frac{\partial}{\partial u_{is}}+u_{ij} \frac{\partial}{\partial u_{ii}},
 $$
 $$
  P_{ij}=2\sum_s u_{is}u_{js}\frac{\partial}{\partial u_{ss}}+\sum_{s\ne k}u_{is}u_{jk} \frac{\partial}{\partial u_{sk}};
 $$
 notice that $X_{ij}=X_{ji}$ and $P_{ij}= P_{ji}$, while  $L_{ij}\ne L_{ji}$.  Thus, we have
 $n(n+1)/2$ operators $X_{ij}$,  $n^2$ operators $L_{ij}$ and  $n(n+1)/2$ operators $P_{ij}$.
 Altogether, they generate a Lie algebra $Sp(2n)$ of dimension $n(2n+1)$. We will work with the
 affine Pl\"ucker embedding of the Lagrangian Grassmannian $\Lambda$ in the projective
 space $P^{N-1}$ specified by the position vector
 $$
 {\bf r}=(M_1, M_2, ..., M_n);
 $$
 here  $M_l$ denotes a  collection of all $l\times l$ minors of $U$, and $N=C_{2n}^n-C_{2n}^{n+2}$ is the
 total number of minors of $U$ (equivalently, $N$ is the dimension of the vector space of effective $n$-forms).
Thus, in the case $n=3$ one has
 $$
 {\rm dim} M_1=6, \ {\rm dim} M_2=6, \ {\rm dim} M_3=1, ~~~~ \Lambda^6 \subset P^{13}.
 $$
 Similarly, for $n=4$ one has
  $$
 {\rm dim} M_1=10, \ {\rm dim} M_2=20, \ {\rm dim} M_3=10, \ {\rm dim} M_4=1, ~~~~ \Lambda^{10} \subset P^{41}.
 $$
With each point of the Pl\"ucker embedding of $\Lambda$ one can
associate a sequence of osculating subspaces
$$
O_1\subset O_2\subset O_3\subset ... \subset O_{n}=P^{N-1}
$$
where $O_1$ is the tangent space, $O_2$ is the second osculating subspace, etc.
At the point corresponding to $u_{ij}=0$ the subspace $O_1$ is spanned by
vectors from  $M_1$, the subspace $O_2$ is spanned by vectors from $M_1\cup
M_2$, etc.
 A symplectic Monge-Amp\'ere equation is naturally identified with a hyperplane section of $\Lambda$.
Below we present three equivalent characterisations of
linearisable equations.

\medskip

 \noindent {\bf Theorem 2.} {\it For a non-degenerate symplectic Monge-Amp\'ere equation in $n$ dimensions the
 following conditions are equivalent:

 \noindent (1) The equation is linearisable by a transformation from $Sp(2n)$.

 \noindent (2) The equation is invariant under  an $n^2$-dimensional subalgebra of $Sp(2n)$.

 \noindent (3) The equation corresponds to a hyperplane which contains an osculating space  $O_{n-2}$.
}

\medskip

 \centerline{\bf Proof:}

 \medskip

 \noindent {\bf equivalence $(1)\Longleftrightarrow (2)$}: since all $n$-dimensional linear equations
 are equivalent under the (complexified) action of $Sp(2n)$, let us consider  the Laplace equation,
 $$
 u_{11}+\dots +u_{nn}=0.
 $$
A direct calculation shows that it is invariant under the
$n^2$-dimensional subalgebra of $Sp(2n)$  generated by the
operators $X_{ij}$ ($i\ne j$), $X_{ii}-X_{jj}$, $L_{ij}-L_{ji}$
($i\ne j$), and $\sum L_{ii}$. Since this property is manifestly
$Sp(2n)$-invariant, any non-degenerate linearisable equation is
invariant under  a subalgebra of $Sp(2n)$ of dimension $n^2$. The
converse is true for any (not necessarily Monge-Amp\'ere)
non-degenerate equation of the form $F(u_{ij})=0$. Indeed, let $G$
be a symmetry group of such equation. We can always assume that
the point $0$, specified by $u_{ij}=0$, belongs to the
hypersurface in the Lagrangian Grassmannian $\Lambda$
corresponding to our equation. Let $G_0$ be the stabilizer of this
point in $G$. Note that $\dim G - \dim G_0 \le \dim
\Lambda-1=n(n+1)/2-1$, as $G$ takes the equation to itself. The
stabilizer of the point $0\in \Lambda$ under the action of
$Sp(2n)$ has the form:
\[
P = \left\{\left.\begin{pmatrix} A & AB \\ 0 & {}^t\!A^{-1} \end{pmatrix}
\,\right|\, A\in GL(n), B={}^t\!B\right\}.
\]
Since the equation is non-degenerate, we can always bring it to the
canonical form:
\begin{equation}\label{cf}
F(u_{ij}) = u_{11} + \dots + u_{nn} + o(u_{ij}) = 0.
\end{equation}
This form (together with the point $0$) is stabilized by the elements of $P$
given by the condition $A\in CO(n), \ B=0$. Hence,  $\dim G_0\leq \dim CO(n)$ so that $\dim G
\le \dim \Lambda - 1 + \dim CO(n) = n^2$. The equality holds only if
$ G_0= CO(n)$. But the subgroup of scalar matrices in $CO(n)$ acts by non-trivial
scalings on the terms of the order 2 and higher in \eqref{cf}. Hence, if
$ G_0= CO(n)$, all higher order terms vanish identically, and our equation coincides
with the Laplace equation.

\medskip

\noindent {\bf equivalence $(1)\Longleftrightarrow (3)$}: Let us
fix a point on $\Lambda$ corresponding to the parameter values
$u_{ij}=0$ (since $Sp(2n)$ acts transitively on $\Lambda$, all
points are equivalent). The osculating  space $O_{n-2}$ at this
point is  specified by the equations $M_{n-1}=0, \ M_n=0$. Any
hyperplane containing this osculating space corresponds to an
equation which is a linear combination of  minors of $U$ of the
orders $n-1$ and $n$. Any such equation  linearises under  the
Legendre transformation. This finishes the proof.

%\medskip \noindent {\bf Remark 1.} Any hyperplane $H$ containing an osculating subspace $O_{n-2}$ is tangential to the Lagrangian Grassmannian $\Lambda$ along a submanifold of codimension 6.  In coordinates $u_{ij}$ and an osculating subspace at the origin, it  coincides with all $n\times n$ symmetric matrices of rank $n-3$.

\medskip

\noindent {\bf Remark.} Linearisable equations of Monge-Amp\'ere
type are of interest in their own, for instance, the equation
$$
u_{tt}(1+u_{xx}+u_{yy})-u_{xt}^2-u_{yt}^2=\epsilon
$$
governs the evolution  of `K\"ahler potentials'  \cite{Donaldson}. It corresponds
to the Lagrangian
$$
\int \left[u_t^2(1+u_{xx}+u_{yy})+2\epsilon u\right] dxdydt,
$$
and linearises under the partial Legendre transformation, $\tilde x=x,\  \tilde
y=y,\  \tilde t=u_t,\  \tilde u=tu_t-u, \ \tilde u_{\tilde x}=-u_x, \  \tilde
u_{\tilde y}=-u_y,\  \tilde u_{\tilde t}=t$, taking the form
$$
1-\tilde u_{\tilde x\tilde x}- \tilde u_{\tilde y\tilde y }=\epsilon \tilde
u_{\tilde t\tilde t}.
$$
Notice that all symplectic Monge-Amp\'ere equations are
Lagrangian. This is a general fact which does not require the
linearisability/integrability.

\section{Classification of integrable  equations}

In this section we classify integrable four-dimensional symplectic
Monge-Amp\'ere equations up to the action of the equivalence group
$Sp(8)$. In Sect.~3.1 we begin with a geometric discussion which
establishes a correspondence between  integrable equations and
`maximally singular' hyperplane sections of the Lagrangian
Grassmannian $\Lambda^{10}\subset P^{41}$.  This leads to the main
classification result proved in Sect.~3.2. Lax pairs, symmetry
algebras and singular varieties associated with integrable
equations arising in the classification are discussed in Sect.~3.3
-- 3.5.

\subsection{Geometric preliminaries: proof of Theorem \ref{thm3}}

Let $v\in V$ be a vector in the $8$-dimensional linear symplectic
space $V$. Lagrangian subspaces passing through $v$ form a
sub-Grassmannian $\Lambda^6\subset \Lambda^{10}\subset P^{41}$.
Thus, the Lagrangian Grassmannian $\Lambda^{10}$ is foliated by a
7-parameter family of $\Lambda^{6}$. Each $\Lambda^6$ lies in some
$P^{13}\subset P^{41}$. Given a four-dimensional symplectic
Monge-Amp\'ere equation (\ref{M}), let $\pi$ be the corresponding
hyperplane in $P^{41}$. The condition that all traveling wave
reductions of the four-dimensional equation are linearisable is
equivalent to the requirement that the corresponding hyperplane
$\pi$ is tangential to all $\Lambda^6\subset \Lambda^{10}$. This
can only happen if either of the following conditions is met:

\noindent (1) The $7$-parameter family of sub-Grassmannians
$\Lambda^6$ possesses an envelope $E\subset \Lambda^{10}$ which is
contained in $\pi$;

\noindent (2) The hyperplane $\pi$ is tangential to $\Lambda^{10}$ along a
$4$-dimensional subvariety $X^4\subset \Lambda^{10}$ which meets each
$\Lambda^6$ (we point out that no hyperplane can be tangential to
$\Lambda^{10}$ along a subvariety of more than four dimensions).

\noindent The first case can be ruled out straightaway: the
$7$-parameter family of sub-Grassmannians $\Lambda^6$ is `too big'
to possess a non-trivial envelope (e.g., a one-parameter family of
lines in the plane always possesses an envelope, however, the full
two-parameter family of lines doesn't). Thus, we are left with the
second possibility, which establishes the necessity part of
Theorem~\ref{thm3}. The sufficiency will follow from the
classification of all non-degenerate equations with a
four-dimensional variety of tangency provided in Sect. 3.2: the
requirement of the existence of such $X^4$ proves to be very
strong indeed, and leads to a finite list of examples all of which
turn out to be integrable.

In the following we will call the variety of tangency $X^4$
\emph{the singular variety} of the symplectic Monge-Amp\'ere
equation defined by the hyperplane section $\pi$. Indeed, $X^4$
coincides with  the set of singular points of the intersection of
$\pi$ with $\Lambda^{10}$.

The condition that the singular variety $X^4$ meets \emph{each}
$\Lambda^6$ is essential for the integrability. In more detail,
consider the incidence variety $Y\subset \Lambda^{10} \times PV$,
\[
 Y = \{ (l, \langle v \rangle) \mid v\in l \},
\]
along with the two projections $\pi_1\colon Y\to \Lambda^{10}$ and
$\pi_2\colon Y\to PV$. Then $\pi_1^{-1}(X^4)$ is a 7-dimensional
closed algebraic variety, and the projection $\pi_2(\pi_1^{-1}(X^4))$ is a
closed algebraic subvariety in $PV$. The condition that $X^4$ meets all
sub-Grassmannians $\Lambda^6$ is equivalent to $\pi_2(\pi_1^{-1}(X^4))=PV$.
This allows us to check this condition using affine charts on $\Lambda^{10}$ which
have non-empty intersections with $X^4$.

\medskip
\noindent\textbf{Example 1.}  Let us describe the singular variety $x^4$
for the non-degenerate linear equation
 $u_{12}=u_{34}$. As shown above, the corresponding
hyperplane has a contact of order $2$ with the Lagrangian
Grassmannian $\Lambda^{10}$ at the Lagrangian subspace $l_\infty$
given by $x^i=0$. There is also a natural conformal structure on
$l_\infty$ associated with this equation. Explicitly, it is given
by $du_1du_2-du_3du_4=0$. A direct computation in various affine
charts of $\Lambda^{10}$ shows that the set of points where the
corresponding hyperplane is tangential to $\Lambda^{10}$ consists
of all Lagrangian subspaces $l$ such that one of the two following
cases holds:
\begin{itemize}
\item $l$ has a 3-dimensional intersection with $l_{\infty}$;
\item $l$ has a 2-dimensional intersection with $l_{\infty}$
    which is isotropic with respect to the above conformal
    structure.
\end{itemize}
Altogether, these points form a 4-dimensional closed algebraic
variety with three irreducible components: one for subspaces with
a 3-dimensional intersection, and  two for subspaces with
2-dimensional isotropic intersections.

It is easy to see that the first component meets all $\Lambda^6$,
that is for any $v\in V$ there exists a Lagrangian subspace
$l\subset V$ such that $v\in l$ and $\dim l\cap l_{\infty}=3$.
Moreover, if $v$ does not lie in $l_{\infty}$, such Lagrangian
subspace $l$ is unique.

\medskip

\noindent\textbf{Example 2.}  Consider the equation $\Hess u=1$, which is
known to be non-integrable. Let us show that its singular variery $X^4$ misses most of the $\Lambda^6$.
Partial Legendre transform brings it to $u_{11}u_{22}-u_{12}^2 = u_{33}u_{44}-u_{34}^2$. The set of
points where the corresponding hyperplane is tangential to $\Lambda^{10}$
is given by the equations $u_{11}=u_{22}=u_{12}=0$ and
$u_{33}=u_{44}=u_{34}=0$. In other words, the intersection $X^4_0$
of $X^4$ with an affine chart of $\Lambda^{10}$ defined by
subspaces $x^i=0$ and $u_j=0$ is represented as the following set
of symmetric matrices,
\[
\begin{pmatrix} 0 & 0 & u_{13} & u_{14} \\
0 & 0 & u_{23} & u_{24} \\
u_{13} & u_{23} & 0 & 0 \\
u_{14} & u_{24} & 0 & 0
\end{pmatrix}.
\]
It is easy to check that $\pi_2(\pi_1^{-1}(X^4_0))$ lies in a quadric in $PV$
and, hence, $\pi_2(\pi_1^{-1}(X^4_0))\ne PV$. This
confirms  that the equation $\Hess u = 1$ is not integrable, even though the
corresponding hyperplane is tangential to $\Lambda^{10}$ along a 4-dimensional
variety.

\medskip

As we have seen in Example~1, the variety of tangency $X^4$  may
consist of several irreducible components with different geometric
properties. However, if the equation under study is non-degenerate
and non-linearisable, any two generic Lagrangian subspaces from
$X^4$ must have trivial intersections. This follows from the two
Lemmas which analyse the cases when  generic pairs of Lagrangian
subspaces of $X^4$ have two- or one-dimensional intersections.

\medskip

\noindent {\bf Lemma 1}. Let $X^4$ be a four-dimensional
subvariety of $\Lambda^{10}$  such that any two generic Lagrangian
subspaces from $X^4$ have two-dimensional intersections. Then
there exists a fixed Lagrangian subspace $L$ which has
three-dimensional intersections with all subspaces from $X^4$.
This situation corresponds to linearisable equations.

\medskip

\centerline{\bf Proof:} Let $l_1, l_2, l_3$ be any three generic
subspaces from $X^4$ such that ${\rm dim} (l_i\cap l_j)=2$. Then
${\rm dim} (l_1\cap l_2\cap l_3)=1$, and we can introduce the
Lagrangian subspace $l={\rm span} \{l_1\cap  l_2, l_1\cap  l_3,
l_2\cap l_3\}$. Any other generic Lagrangian subspace from $X^4$
will have one-dimensional intersections with $l_1\cap l_2, \
l_1\cap l_3$ and $l_2\cap l_3$. Thus, it will have a
three-dimensional intersection with $l$.
%(unless it passes through $l_1\cap l_2\cap l_3$, which is forbidden due to the non-degeneracy).
Any hyperplane $\pi$ which is tangential to $\Lambda^{10}$ along
$X^4$ will automatically contain the osculating subspace to
$\Lambda^{10}$ at $l$, so that the corresponding equation will be
linearisable.

\medskip

\noindent {\bf Lemma 2}. There exists no non-degenerate symplectic
Monge-Amp\'ere equation for which any two generic Lagrangian
subspaces in the singular variety $X^4$ have one-dimensional
intersections.

\medskip

\centerline{\bf Proof:} Let $l_1$ and $l_2$ be two Lagrangian
subspaces from $X^4$ such that ${\rm dim} (l_1\cap l_2)=1$. Let
$v$ be a vector in the intersection. Consider the Lagrangian
Grassmannian $\Lambda^6\subset P^{13}$ generated by $v$. Then the
travelling wave reduction in the direction of $v$ produces a
three-dimensional linearisable Monge-Amp\'ere equation such that
the corresponding hyperplane in $P^{13}$ is tangential to
$\Lambda^6$ at two distinct points, $\tilde l_1=l_1/\langle v\rangle$ and
$\tilde l_2= l_2 /\langle v\rangle$. The points $\tilde l_1$ and $\tilde l_2$
correspond to three-dimensional Lagrangian subspaces in $V^6$
which do not intersect by construction. This situation is,
however, not possible for linearisable equations in three
dimensions (even degenerate). Thus, the travelling wave reduction
in the direction of $v$ must produce the `zero' equation, that is,
 $\Lambda^6\subset \pi\cap \Lambda^{10}$ where $\pi$ is the
 hyperplane corresponding to our Monge-Amp\'ere equation. When $l_1$
 and $l_2$ are allowed to vary within $X^4$, one obtains at least a
 four-parameter variety of pairwise intersections $v$
 (as explained below). Thus, we obtain a four-parameter
 family of  $\Lambda^6$  which belong to
 the nine-dimensional hyperplane section $\pi\cap \Lambda^{10}$. This,
 however, is not possible since the union of any four-parameter
 family of $\Lambda^6$ must be essentially ten-dimensional.

 It remains to prove that we have at least a four-dimensional
 variety of  pairwise intersections $v$. Let us fix $l_1$
 and consider the set of all intersections $Y(l_1)=\{l_1\cap l_2, \ l_2\in X^4\}$. Let $s$ be the dimension of $Y(l_1)$.
 Fixing $v$ in $Y(l_1)$ we obtain a $(4-s)$-parametric family of
 Lagrangian subspaces $l$ in $X^4$ which contain $v$. Each of them
 contains an s-parameter subvariety $Y(l)$. Taking a union of $Y(l)$ over all $l$ passing through $v$, one
 obtains a four-parameter family of pairwise intersections $v$.
 This finishes the proof of Lemma 2.

\medskip

Let us now return to the general case when any two generic
Lagrangian subspaces of $X^4$ have trivial intersections. Take any
two such planes, say $l_1$ and $l_2$, and introduce canonical
coordinates $x^i, \ u_i$ in $V$ such that $l_1, l_2$ become the
coordinate planes $x^i=0$ and $u_i=0$, respectively. Suppose now
that a symplectic Monge-Amp\'ere equation corresponds to a
hyperplane $\pi$ which is tangential to the Lagrangian
Grassmannian $\Lambda^{10}$ at the points corresponding to $l_1$
and $l_2$. In coordinates $x^i, u_i$, any such equation takes a
`purely quadratic' form generated by a linear combination  of
$2\times 2$ minors of the matrix $U$ only. Notice that such form
is invariant under $GL(4)\subset Sp(8)$. Let us choose a third
Lagrangian plane $l_3$ which does not intersect $l_1$ and $l_2$,
and bring it to a normal form using the remaining $GL(4)$-freedom.
The requirement that $\pi$ is tangential to $\Lambda^{10}$ at the
point corresponding to $l_3$ imposes further constraints on the
quadratic part of the equation. Further details of the
classification are provided in the next Section.

\subsection{Classification results: proof of Theorem \ref{thm4}}

\medskip

\noindent{\bf Theorem 4.} {\it Over the field of complex numbers any integrable
non-degenerate symplectic  Monge-Amp\'ere equation is $Sp(8)$-equivalent to one
of the following  normal forms:
\begin{enumerate}
\item $u_{11}-u_{22}-u_{33}-u_{44}=0$ (linear wave equation);
\item $u_{13}+u_{24}+u_{11}u_{22}-u_{12}^2=0$ (second heavenly equation);
\item $u_{13} = u_{12}u_{44}-u_{14}u_{24}$ (modified heavenly equation);
\item $u_{13}u_{24}-u_{14}u_{23}=1$ (first heavenly equation).
\item $u_{11}+u_{22}+u_{13}u_{24}-u_{14}u_{23}=0$ (Husain equation);
\item $\alpha u_{12}u_{34} + \beta u_{13}u_{24}+\gamma u_{14}u_{23} = 0$ (general
heavenly equation),  $\alpha+\beta+\gamma=0$.
\end{enumerate}
}

\medskip

\centerline{\bf Proof:} As shown above,  we can always assume that
the Monge-Amp\'ere equation contains only quadratic terms, or, in
other words, the corresponding hyperplane $\pi\subset P^{41}$ is
tangential to the Lagrangian Grassmannian $\Lambda^{10}$ at the
points $l_1 =\infty$ (specified by $x^i=0$) and $l_2=0$ (specified
by $u_i=0$). Any other Lagrangian plane $l_3$ which has trivial
intersections with $l_1$ and $l_2$ is defined by a non-degenerate
symmetric $4\times4$ matrix $u_{ij}$. It can be viewed up to the
natural action of the subgroup $GL(4)\subset Sp(8)$ which
stabilizes $l_1$ and $l_2$. Over the field of complex numbers all
such symmetric matrices are equivalent to each other. Fix $l_3$
such that $u_{14}=1$, $u_{23}=-1$ and all other $u_{ij}$ equal
$0$. Assume that, in addition to $l_1$ and $l_2$, our hyperplane
$\pi$ is also tangential to $\Lambda^{10}$
 at the point $l_3$.

The set of all quadratic Monge-Amp\'ere equations such that the corresponding
hyperplanes are tangential to the Lagrangian Grassmannian $\Lambda^{10}$ at the
points $l_1,l_2,l_3$ forms a $10$-dimensional vector space spanned by the
following expressions:
\begin{align*}
E_0 &= u_{11}u_{22}-u_{12}^2,\\
E_1 &= 1/2(u_{11}u_{24}-u_{12}u_{14}+u_{22}u_{13}-u_{12}u_{23}),\\
E_2 &= 1/6(u_{11}u_{44}-u_{14}^2+u_{22}u_{33}-u_{23}^2)+1/3(2u_{13}u_{24}-u_{14}u_{23}-u_{12}u_{34}),\\
E_3 &= 1/2(u_{33}u_{24}-u_{23}u_{34}+u_{44}u_{13}-u_{14}u_{34}),\\
E_4 &= u_{33}u_{44}-u_{34}^2;
\end{align*}
and
\begin{align*}
F_0 &= u_{11}u_{33}-u_{13}^2,\\
F_1 &= 1/2(u_{11}u_{34}-u_{13}u_{14}+u_{33}u_{12}-u_{13}u_{23}),\\
F_2 &= 1/6(u_{11}u_{44}-u_{14}^2+u_{22}u_{33}-u_{23}^2)+1/3(2u_{12}u_{34}-u_{14}u_{23}-u_{13}u_{24}),\\
F_3 &= 1/2(u_{22}u_{34}-u_{23}u_{24}+u_{44}u_{12}-u_{14}u_{24}),\\
F_4 &= u_{22}u_{44}-u_{24}^2.
\end{align*}
Denote  by $E$ and $F$ the subspaces spanned by $E_i$ and $F_j$,
respectively. Note that permutations of indices $(1,2)$ and
$(3,4)$ leave both subspaces invariant, while permutations $(1,4)$
and $(2,3)$ interchange them. Both subspaces are also invariant
under the natural action of $SO(4,\bbbc)$ given as the set of all
linear symplectic transformations stabilizing  $l_1$, $l_2$ and
$l_3$. Moreover, it is possible to show that under the
identification of $SO(4,\bbbc)$ with $(SL(2,\bbbc)\times
SL(2,\bbbc))/{(\pm E_2, \pm E_2)}$, each of the copies of $SL(2,
\bbbc)$ preserves the decomposition $E\oplus F$, acting
irreducibly on one summand and trivially on the other.

The standard model for an irreducible 5-dimensional $SL(2,\bbbc)$-action is
given as the following action  on degree 4 polynomials:
\[
\begin{pmatrix}
a & b \\ c & d
\end{pmatrix}\ p(t) = (ct+d)^4 p\left(\frac{at+b}{ct+d}\right).
\]
The vectors $E_i$ and $F_j$  are chosen in such a way that the
correspondence $E_i\leftrightarrow t^i$ and $F_i\leftrightarrow
t^i$, $i=0,\dots,4$,  establishes the isomorphisms of
$SL(2,\bbbc)$-actions on subspaces $E$ and $F$ with the standard
model. Polynomials of degree 4 are easily classified with respect
to the $SL(2,\bbbc)$ action and the multiplication by a non-zero
scalar. We have the following set of canonical representatives
according to the roots of the polynomial:
\begin{enumerate}
\item[(0)] $0$;
\item[(1)] $1$ (all four roots of $p(t)$ coincide and are equal to $\infty$);
\item[(2)] $t$ ($p(t)$ has one triple root $\infty$ and one simple root $0$);
\item[(3)] $t^2$ ($p(t)$ has two double roots $0$ and $\infty$);
\item[(4)] $t^2-1$ ($p(t)$ has one double root $\infty$ and two simple roots $\pm1$);
\item[(5)] $(t^2-1)(t-a)$, $a\ne\pm 1$ ($p(t)$ has four simple roots $\infty$, $\pm1$, $a$).
\end{enumerate}
This classification gives us an easy way to classify vectors in
$E\oplus F$ viewed up to the action of $SO(4,\bbbc)$, permutation
of $E$ and $F$ and the multiplication by a non-zero scalar.
Namely, decompose any non-zero vector $v$ in $E\oplus F$ into the
difference $v_e-v_f$, $v_e\in E,v_f\in F$ and match each of these
two vectors with the corresponding polynomial in the standard
model. This gives a pair of polynomials $(p(t),q(t))$, which fully
describe the vector $v$. Next, we classify this pair up to
independent actions of $SL(2,\bbbc)$ on each of them, permutations
$p\leftrightarrow q$ and a multiplication by a non-zero constant.
We easily get the set of normal forms $v_{i,j}=(p_i(t),
cq_j(t))$,\ $i,j=0,\dots,5$,\ $i\ge j$,\ $c\ne 0$, where $p_i$ and
$q_j$ are polynomials listed above under the $i$-th and $j$-th
items. We can assume  $c\ne 0$ since the case $c=0$ is included in
the subcase $q(t)=q_0(t)=0$. Next, we consider each of these
normal forms case-by-case and check for which of the corresponding
Monge-Amp\'ere equations the variety of  points where the
corresponding hyperplane $\pi$ is tangential to the Lagrangian
Grassmannian $\Lambda^{10}$ has dimension $4$. The complete list
of all such cases is given below:
\begin{enumerate}
\item $p(t)=q(t)=(t^2-1)(t-a)$,~ $a\ne \pm 1$;
\item $p(t)=q(t)=t^2-1$;
\item $p(t)=t^2-1, \ q(t)=t^2$;
\item $p(t)=q(t)=t^2$;
\item $p(t)=q(t)=t$;
\item $p(t)=t, \ q(t)=1$;
\item $p(t)=q(t)=1$;
\item $p(t)=t^3-t,\ q(t)=0$;
\item $p(t)=t,\ q(t)=0$;
\item $p(t)=1,\ q(t)=0$.
\end{enumerate}
The cases 4, 7 and 10 correspond to the degenerate equations
$u_{13}u_{24}-u_{12}u_{34}=0$,
$u_{11}u_{22}-u_{12}^2=u_{11}u_{33}-u_{13}^2$ and
$u_{11}u_{22}-u_{12}^2=0$, respectively. In fact, the last two
equations are equivalent to the degenerate linear equations
$u_{22}=u_{33}$ and $u_{22}=0$, respectively, via the partial
Legendre transform of $(x^1,u_1)$: $(x^1,u_1)\mapsto (-u_1,x^1)$.

The case 9 is given by
\[
u_{11}u_{24}-u_{12}u_{14}+u_{22}u_{13}-u_{12}u_{23} = 0,
\]
and is equivalent to a non-degenerate linear equation (and, hence, to the
linear wave equation) via the partial Legendre transforms of $(x^1,u_1)$ and
$(x^2,u_2)$.

Similarly, the case~8 is equivalent to $\Hess u = 1$. Indeed,
the polynomial $t^3-t$ is equivalent to $t^4-1$ under the action
of $SL(2,\bbbc)$ since in both cases all 4 roots are distinct and
their cross-ratio is equal to $-1$. The pair of polynomials
$p(t)=t^4-1$, $q(t)=0$ corresponds to the equation
$u_{11}u_{22}-u_{12}^2=u_{33}u_{44}-u_{34}^2$ which can be
brought to $\Hess u = 1$ via the partial Legendre transform of
$(x^1,u_1)$ and $(x^2,u_2)$. As we have seen in Example~2,
this equation is not integrable.

Other cases can be considered in a similar manner. It appears that
they correspond  to the remaining non-linear integrable
Monge-Amp\'ere equations listed in Theorem~\ref{thm4}. Namely,
equation~5 can be brought to the modified heavenly equation via
the change of variables $(x^1,x^2,x^3,x^4)\mapsto
(x^1,x^2+x^3,x^2-x^3,x^4)$, partial Legendre transform of
$(x^1,u_1)$ and a permutation of indices.

Equations from cases~2 and~3 are brought to the Husain  and the first
heavenly equations, respectively, via the partial Legendre transform of
$(x^1,u_1)$.

The case 6 is brought to the second heavenly equation via the partial
Legendre transforms of $(x^1,u_1)$ and $(x^2,u_2)$, and a permutation of
indices.

Finally, to establish the correspondence between equation~1 and the general
heavenly equation one needs to transform the canonical form $(t^2-1)(t-a)$ to
$t^4+2bt^2+1$ under the action of $SL(2,\bbbc)$ and then use the following
change of variables: $(x^1,x^2,x^3,x^4)\mapsto
(x^1+x^4,x^2+x^3,x^2-x^3,x^1-x^4)$.

The obtained correspondence is summarized in the following table:
\begin{center}
\begin{tabular}{|c|c|c|l|}
\hline
Case & $p(t)$ & $q(t)$ & Equation \\
\hline 1 & $(t^2-1)(t-a)$ & $(t^2-1)(t-a)$ & \emph{general heavenly} \\
\hline 2 & $t^2-1$ & $t^2-1$ & \emph{Husain equation} \\
\hline 3 & $t^2-1$ & $t^2$ & \emph{first heavenly} \\
\hline 4 & $t^2$ & $t^2$ & degenerate equation \\
\hline 5 & $t^2$ & $t$ & \emph{modified heavenly} \\
\hline 6 & $t$ & $1$ & \emph{second heavenly} \\
\hline 7 & $1$ & $1$ & degenerate equation \\
\hline 8 & $t(t^2-1)$ & 0 & $\Hess u = 1$ (non-integrable) \\
\hline 9 & $t$ & $0$ & \emph{linear wave equation} \\
\hline 10 & 1 & 0 & degenerate equation \\
\hline
\end{tabular}
\end{center}
It remains to point out that all equations appearing in Theorem
\ref{thm4} are indeed integrable: they possess the required number
of hydrodynamic reductions (for the first and second heavenly
equations this has been demonstrated in \cite{Fer2, Fer3}), along
with the Lax pairs depending on an auxiliary spectral parameter.
This completes the proof of Theorem~\ref{thm4}.

\medskip
\noindent {\bf Remark.} Although all equations listed in Theorem \ref{thm4} are
not $Sp(8)$-equivalent (in fact, not even contact equivalent), there may exist
more complicated connections among them. Below we describe the construction of
\cite{Plebanski} which links the first and second heavenly equations. Let us
take the second heavenly equation in the form $\theta_{1\tilde
1}+\theta_{2\tilde 2}+\theta_{\tilde1 \tilde1}\theta_{\tilde 2
\tilde2}-\theta_{\tilde 1\tilde 2}^2=0$, here $\theta$ is a function of the
four independent variables $x^1,  x^2, {\tilde x^1}, {\tilde x^2}$. Introducing
the two-form
$$
\Omega=(d{\tilde x^1}-\theta_{\tilde 2\tilde 2}dx^1+\theta_{\tilde 1\tilde
2}dx^2)\wedge (d{\tilde x^2}+\theta_{\tilde 1\tilde 2}dx^1-\theta_{\tilde
1\tilde 1}dx^2),
$$
one can verify the relations $d\Omega=0, \ \Omega \wedge \Omega=0$. Thus, by
Darboux's theorem, there exist coordinates $x^3, x^4$ such that $\Omega
=dx^3\wedge dx^4$. This implies the expansions
$$
d{\tilde x^1}-\theta_{\tilde 2\tilde 2}dx^1+\theta_{\tilde 1\tilde
2}dx^2=u_{13}dx^3+u_{14}dx^4, ~~~~ d{\tilde x^2}+\theta_{\tilde 1\tilde
2}dx^1-\theta_{\tilde 1\tilde 1}dx^2=u_{23}dx^3+u_{24}dx^4
$$
where $u_{ij}$  are so far arbitrary coefficients satisfying the relation
$u_{13}u_{24}-u_{14}u_{23}=1$. Rewriting these relations in the form
$$
d{\tilde x^1}=\theta_{\tilde 2\tilde 2}dx^1-\theta_{\tilde 1\tilde
2}dx^2+u_{13}dx^3+u_{14}dx^4, ~~~~ d{\tilde x^2}=-\theta_{\tilde 1\tilde
2}dx^1+\theta_{\tilde 1\tilde 1}dx^2+u_{23}dx^3+u_{24}dx^4,
$$
one can see that, by virtue of the relation $\partial {\tilde x^1}/\partial
x^2=\partial {\tilde x^2}/\partial x^1$, there exists a function $u(x^1, x^2,
x^3, x^4)$ such that ${\tilde x^1}=\partial u/\partial x^1, \ {\tilde
x^2}=\partial u/\partial x^2$. In this case coefficients $u_{ij}$ become second
order partial derivatives of $u$, and the relation
$u_{13}u_{24}-u_{14}u_{23}=1$ becomes the second heavenly equation. The passage
from $x^1,  x^2, {\tilde x^1}, {\tilde x^2}$ to the new independent variables
$x^1,  x^2, x^3, x^4$ can be viewed as a multi-dimensional analogue of
reciprocal transformations familiar from the $(1+1)$-dimensional theory.  This
link is highly nonlocal, and in any case not a contact equivalence. The only
`local' relations one can write here are the following:
$$
{\tilde x^1}=u_1, \  {\tilde x^2}=u_2, \ \theta_{\tilde 2\tilde 2}=u_{11}, \
\theta_{\tilde 1\tilde 2}=-u_{12}, \ \theta_{\tilde 1\tilde 1}=u_{22}.
 $$
One can expect similar non-local links between other equations listed in
Theorem \ref{thm4}.

\subsection{Lax pairs}

In this section we present Lax pairs for integrable equations from Theorem 4.
Some of them readily follow from the Lax pair for the six-dimensional version
of the second heavenly equation \cite{Takasaki1, Przanovski},
$$
u_{15}+u_{26}+u_{13}u_{24}-u_{14}u_{23}=0,
$$
which is given by a pair of vectors fields
$$
X_1=\partial_6+u_{13}\partial_4-u_{14}\partial_3+\lambda \partial_1, ~~~
X_2=\partial_5-u_{23}\partial_4+u_{24}\partial_3-\lambda \partial_2,
$$
$\lambda=$const, such that  $[X_1, X_2]=0$ modulo the equation.
Straightforward dimensional reductions provide Lax pairs for the
first heavenly, second heavenly and Husain equations (we point out
that the general heavenly equation {\it does not} arise as a
travelling wave reduction of this six-dimensional equation).

\noindent {\bf Second heavenly equation}
$u_{13}+u_{24}+u_{11}u_{22}-u_{12}^2=0$:
$$
X_1=\partial_4+u_{11}\partial_2-u_{12}\partial_1+\lambda \partial_1, ~~~
X_2=\partial_3-u_{12}\partial_2+u_{22}\partial_1-\lambda \partial_2.
$$
\noindent {\bf Modified heavenly equation} $u_{13} =
u_{12}u_{44}-u_{14}u_{24}$:
$$
X_1=u_{14}\partial_2-u_{12}\partial_4+\lambda \partial_1, ~~~
X_2=-\partial_3+u_{44}\partial_2-u_{24}\partial_4+\lambda \partial_4.
$$
\noindent {\bf First heavenly equation} $u_{13}u_{24}-u_{14}u_{23}=1$:
$$
X_1=u_{13}\partial_4-u_{14}\partial_3+\lambda \partial_1, ~~~
X_2=-u_{23}\partial_4+u_{24}\partial_3-\lambda \partial_2.
$$
%\noindent {\bf Grant equation} $u_{11}+u_{13}u_{24}-u_{14}u_{23}=0$:
%$$
%X_1=u_{13}\partial_4-u_{14}\partial_3+\lambda \partial_1, ~~~
%X_2=\partial_1-u_{23}\partial_4+u_{24}\partial_3-\lambda \partial_2.
%$$

\noindent {\bf Husain equation} $u_{11}+u_{22}+u_{13}u_{24}-u_{14}u_{23}=0$:
$$
X_1=\partial_2+u_{13}\partial_4-u_{14}\partial_3+\lambda \partial_1, ~~~
X_2=\partial_1-u_{23}\partial_4+u_{24}\partial_3-\lambda \partial_2.
$$

\noindent {\bf General heavenly equation} $\alpha u_{12}u_{34} + \beta
u_{13}u_{24}+\gamma u_{14}u_{23} = 0$:
$$
X_1=u_{34}\partial_1-u_{13}\partial_4+\gamma \lambda
(u_{34}\partial_1-u_{14}\partial_3),  ~~~
X_2=u_{23}\partial_4-u_{34}\partial_2+\beta \lambda
(u_{34}\partial_2-u_{24}\partial_3).
$$
Notice that in the last case the commutator of $X_1, X_2$ is not identically
zero on the solutions, this condition is substituted by  $[X_1, X_2]=0$ mod $X_1, X_2$. Thus, the
2-dimensional distribution spanned by $X_1, X_2$ is integrable. Modifications
of the inverse scattering transform and  the  $\overline\partial$-dressing method
for Lax pairs of this type were developed in \cite{Manakov, Manakov3,
Bogdanov2}.
%We recall that the Grant equation is equivalent to the first heavenly equation under a Legendre transform and, therefore, is not included in the classification.

\subsection{Symmetry algebras}
In this section we describe  symmetry algebras for equations from
Theorem~\ref{thm4}. To be precise, we consider only those
symmetries which belong to the equivalence group $Sp(8)$ or,
equivalently, stabilisers of the corresponding effective forms. It
should be emphasized that the full contact symmetry algebras of
these equations are infinite dimensional. We adopt the notation of
Sect.~2.

\begin{center}
\begin{tabular}{|p{18mm}|c|c|}
\hline Equation & Generators of the symmetry algebra & Dimension \\
 \hline linear wave &  $\begin{smallmatrix}\\ \\
 X_{11}+X_{22},\  X_{11}+X_{33},\  X_{11}+X_{44},\  X_{12},\  X_{13},\
X_{14},\  X_{23},\  X_{24},\  X_{34}, \\
L_{11}+L_{22}+L_{33}+L_{44},\
L_{12}+L_{21},\  L_{13}+L_{31},\  L_{14}+L_{41},\  L_{23}-L_{32},\
L_{24}-L_{42},\  L_{34}-L_{43}
 \end{smallmatrix}$
& $16$ \\
\hline second heavenly &  $\begin{smallmatrix}\\ \\
X_{11}, \  X_{12},\  X_{13}-X_{24}, \ X_{14}, \ X_{22},\  X_{23},\  X_{33}-L_{24}, \
 X_{34}+2L_{23}, \  X_{44}-L_{13}, \\
L_{12}-L_{43}, \ L_{21}-L_{34},\
L_{14}-L_{23}, \ 2L_{11}+L_{22}+L_{44}, \ L_{11}+2L_{22}+L_{33} \end{smallmatrix}$
& $14$ \\
\hline modified heavenly & $\begin{smallmatrix}\\ \\
X_{11}, \ X_{22}, \ X_{23}, \ X_{24}-L_{34}, \ X_{33}, \ X_{34},\
X_{44}+L_{32}, \\
L_{11}, \ L_{22}+L_{33}, \ 2L_{33}+L_{44}, \ L_{24}, \ P_{11},\
P_{44}-2L_{23}
\end{smallmatrix}$
& $13$ \\
\hline first \hbox{heavenly} & $\begin{smallmatrix}\\ \\
X_{11}, \ X_{12}, \ X_{22}, \ X_{33}, \ X_{34}, \ X_{44}, \ L_{12},\\
L_{21}, \ L_{34}, \ L_{43}, \ L_{11}-L_{22}, \ L_{33}-L_{44},\
L_{11}+L_{22}-L_{33}-L_{44}
\end{smallmatrix}$
& $13$ \\
\hline Husain & $\begin{smallmatrix} \\ \\
X_{11}-X_{22}, \ X_{12}, \ X_{33}, \ X_{34}, \ X_{44}, \ L_{11}+L_{22}, \\
L_{12}-L_{21}, \ L_{33}-L_{44}, \ L_{34}, \ L_{43}, \ P_{11}-P_{22}, \ P_{12}\\
\ \\
\end{smallmatrix}$
& $12$ \\
\hline general heavenly & $\begin{smallmatrix} \\ \\
X_{ii}, \ L_{ii}, \ P_{ii}, ~~ i=1,..., 4.
\end{smallmatrix}$
& $12$ \\
\hline
\end{tabular}
\end{center}
The table  shows that  equations from Theorem~\ref{thm4} are not
equivalent  under the action of the symplectic group. Indeed, the
second heavenly equation is the only one with a 14-dimensional symmetry algebra.
Both the modified heavenly and the first heavenly equations have 13-dimensional
symmetry algebras, but only the symmetry algebra of the modified heavenly
equation has an invariant 2-dimensional non-isotropic subspace in $\bbbc^8$.
Similarly, both the general heavenly and the Husain equations have
12-dimensional symmetry algebras, but the symmetry algebra of the general
heavenly equation is reductive, while that of the Husain
equation is not.

\medskip

\noindent {\bf Remark.}
For any effective $n$-form $\omega$ on the $2n$-dimensional symplectic space $V$ with the symplectic 2-form $\Omega$
we can define a bilinear form $B_\omega$ \cite{Lychagin},
\[
B_\omega\colon (X,Y)\mapsto \frac{i_{X}\omega\wedge i_{Y}\omega\wedge
\Omega}{\Omega^n},
\]
which will be symmetric for odd, and skew-symmetric for even values of $n$.
It is known  that $B_\omega$ is proportional to $\Omega$ for all even values of $n$
(we thank Bertrand Banos for pointing this out).
As $\omega$ is only defined up to a scalar multiple, the form $B_\omega$ also
makes sense only up to multiplication by a non-zero constant. Hence, if
$B_\omega = \lambda \Omega$, we can always assume that either  $\lambda = 0$
(i.e., the form $B_\omega$ vanishes identically), or $\lambda=1$ and
$B_\omega=\Omega$. It is clear that the vanishing of $\lambda$ is an invariant
condition for  the form $\omega$. Below we list  the values of $\lambda$ for
integrable Monge-Amp\'ere equations from Theorem~\ref{thm4}:
\begin{center}
\begin{tabular}{|l|c|}
\hline Equation & Coefficient $\lambda$ \\
\hline linear wave & 0 \\
\hline second heavenly & 0 \\
\hline modified heavenly & 0 \\
\hline first heavenly & 1 \\
\hline Husain & 1 \\
\hline general heavenly & 1 \\
\hline
\end{tabular}
\end{center}
For instance, this confirms that the  modified heavenly equation is not contact equivalent to the
first heavenly equation, although both have stabilizers of dimension $13$.

\subsection{Geometry of singular varieties}

We recall that any integrable equation corresponds to a hyperplane
 which is tangential to the Pl\"ucker embedding of the Lagrangian
Grassmannian $\Lambda^{10}$ along a four-dimensional subvariety
$X^4$ which meets all $\Lambda^6\subset \Lambda^{10}$. The second
condition means that the corresponding congruence of Lagrangian
subspaces fills the symplectic space $V^8$. In this section we
provide invariant descriptions of $X^4$ for all examples arising
in the classification.

\medskip

\noindent {\bf Linear equations:} Fix a Lagrangian subspace
$L\subset V^8$. Then $X^4$ consists of all Lagrangian subspaces
which have three-dimensional intersections with $L$. We emphasize
that this constitutes only one irreducible component of the
singular variety. See Example~1 for more details.

\noindent {\bf Second heavenly equation:} Fix a Lagrangian
subspace $L$ and a two-dimensional isotropic plane $l$ such that
$l\subset L\subset V^8$. Let $k$ be the $P^1$ of lines in $l$, and
$K$ the $P^1$ of three-dimensional subspaces in $L$ which contain
$l$. Fix a projective isomorphism between these two $P^1$'s,
$k\leftrightarrow K$. Let $c$ be the union, over $k$,  of
two-dimensional isotropic planes which pass through  $k$ and are
contained in the corresponding $K$. Similarly, let $C$ be the
union, over $K$, of all Lagrangian subspaces which pass through
$K$. Fix a projective isomorphism  $c\leftrightarrow C$ which
covers the isomorphism $k\leftrightarrow K$. Then $X^4$ consists
of all Lagrangian subspaces which contain an isotropic plane from
$c$, and have a three-dimensional intersection with the
corresponding Lagrangian subspace $C$.

\noindent {\bf Modified heavenly equation:} Fix a decomposition
$V^8=V^2+V^6$  into a direct sum of two skew-orthogonal symplectic
spaces. Let $l$ be a two-dimensional isotropic plane in $V^6$, and
$L$ the corresponding co-isotropic four-dimensional subspace,
$l\subset L\subset V^6$. Let $k$ be the $P^1$ of lines in $l$, and
$K$ the $P^1$ of three-dimensional subspaces in $L$ which contain
$l$. Fix a projective isomorphism between these two $P^1$'s,
$k\leftrightarrow K$. Finally, let $F$ be the family of Lagrangian
subspaces in $V^6$ which pass through one of the $k$'s, and have a
two-dimensional intersection with the corresponding $K$. Notice
that $F$ is three-dimensional. Then $X^4$ consists of all
Lagrangian subspaces which have one-dimensional intersections with
$V^2$, and intersect $V^6$ along a Lagrangian subspace from $F$.

\noindent {\bf First heavenly equation:} Fix a decomposition
$V^8=V^4+V^4$  into a direct sum of two skew-orthogonal symplectic
spaces. Fix a two-dimensional Lagrangian plane in each copy of
$V^4$, denote these planes $l_1$ and $l_2$, respectively. Let
$F_1$ be the family of all Lagrangian planes in the first copy of
$V^4$ which have one-dimensional intersections with $l_1$.
Similarly, let $F_2$ be the family of all Lagrangian planes in the
second copy of $V^4$ which have one-dimensional intersections with
$l_2$. Then $X^4$ consists of all Lagrangian subspaces which
intersect each copy of  $V^4$ along two-dimensional Lagrangian
planes from $F_1$ and $F_2$, respectively.

\noindent {\bf Husain equation:} Fix a decomposition $V^8=V^2+
V^2+V^4$  into a direct sum of three skew-orthogonal symplectic
spaces. Fix a two-dimensional Lagrangian plane $l\subset V^4$ and
consider a family $F$ of all Lagrangian planes in $V^4$ which have
one-dimensional intersections with $l$. Then $X^4$ consists of all
Lagrangian subspaces which have one-dimensional intersections with
each copy of  $V^2$, and intersect $V^4$ along a Lagrangian plane
from $F$.

\noindent {\bf General heavenly equation:} Fix a decomposition
$V^8=V^2+ V^2+V^2+V^2$  into a direct sum of four skew-orthogonal
symplectic spaces. Then $X^4$ consists of all Lagrangian subspaces
which have one-dimensional intersections with each copy of $V^2$
(four-secant subspaces).

\medskip

One can show that congruences of Lagrangian subspaces in $V^8$
corresponding to $X^4$ are of class one, that is, a generic vector
in $V^8$ lies in a unique Lagrangian subspace from $X^4$.

\section*{Acknowledgements}
We have benefited from discussions with   Bertrand Banos,  Maciej
Dunajski, James Grant, Valentin Lychagin, Joseph Landsberg, Maxim
Pavlov, Volodya Rubtsov and Artur Sergeev. We thank the LMS for
their financial support of BD to Loughborough making this
collaboration possible. We also thank the ICMS for giving us an
opportunity to spend two weeks  in Edinburgh under the `research
in groups' program where this work has been completed and the
investigation of higher dimensional Monge-Amp\'ere equations was
initiated.

\end{document}